# NONPARAMETRIC ESTIMATION OF AN ADDITIVE MODEL WITH A LINK FUNCTION

By Joel L. Horowitz[1] and Enno Mammen[2]

*Northwestern University and Universität Mannheim*

This paper describes an estimator of the additive components of a nonparametric additive model with a known link function. When the additive components are twice continuously differentiable, the estimator is asymptotically normally distributed with a rate of convergence in probability of $n^{-2/5}$. This is true regardless of the (finite) dimension of the explanatory variable. Thus, in contrast to the existing asymptotically normal estimator, the new estimator has no curse of dimensionality. Moreover, the estimator has an oracle property. The asymptotic distribution of each additive component is the same as it would be if the other components were known with certainty.

**1. Introduction.** This paper is concerned with nonparametric estimation of the functions $m_1, \ldots, m_d$ in the model

$$(1.1) \qquad Y = F[\mu + m_1(X^1) + \cdots + m_d(X^d)] + U,$$

where $X^j (j = 1, \ldots, d)$ is the $j$th component of the random vector $X \in \mathbb{R}^d$ for some finite $d \geq 2$, $F$ is a known function, $\mu$ is an unknown constant, $m_1, \ldots, m_d$ are unknown functions and $U$ is an unobserved random variable satisfying $\mathbf{E}(U|X = x) = 0$ for almost every $x$. Estimation is based on an i.i.d. random sample $\{Y_i, X_i : i = 1, \ldots, n\}$ of $(Y, X)$. We describe an estimator of the additive components $m_1, \ldots, m_d$ that converges in probability pointwise at the rate $n^{-2/5}$ when $F$ and the $m_j$'s are twice continuously differentiable and the second derivative of $F$ is sufficiently smooth. In contrast to previous estimators, only two derivatives are needed regardless of

Received July 2002; revised March 2004.

[1]Supported in part by NSF Grants SES-99-10925 and SES-03-52675, the Alexander von Humboldt Foundation and Deutsche Forschungsgemeinschaft Sonderforschungsbereich 373, "Quantifikation und Simulation Ökonomischer Prozesse."

[2]Supported in part by Deutsche Forschungsgemeinschaft MA 1026/6-2.

the dimension of $X$, so asymptotically there is no curse of dimensionality. Moreover, the estimators derived here have an oracle property. Specifically, the centered, scaled estimator of each additive component is asymptotically normally distributed with the same mean and variance that it would have if the other components were known.

Linton and Härdle (1996) (hereinafter LH) developed an estimator of the additive components of (1.1) that is based on marginal integration. The marginal integration method is discussed in more detail below. The estimator of LH converges at the rate $n^{-2/5}$ and is asymptotically normally distributed, but it requires the $m_j$'s to have an increasing number of derivatives as the dimension of $X$ increases. Thus, it suffers from the curse of dimensionality. Our estimator avoids this problem.

There is a large body of research on estimation of (1.1) when $F$ is the identity function so that $Y = \mu + m_1(X^1) + \cdots + m_d(X^d) + U$. Stone (1985, 1986) showed that $n^{-2/5}$ is the optimal $L_2$ rate of convergence of an estimator of the $m_j$'s when they are twice continuously differentiable. Stone (1994) and Newey (1997) describe spline estimators whose $L_2$ rate of convergence is $n^{-2/5}$, but the pointwise rates of convergence and asymptotic distributions of spline and other series estimators remain unknown. Breiman and Friedman (1985), Buja, Hastie and Tibshirani (1989), Hastie and Tibshirani (1990), Opsomer and Ruppert (1997), Mammen, Linton and Nielsen (1999) and Opsomer (2000) have investigated the properties of backfitting procedures. Mammen, Linton and Nielsen (1999) give conditions under which a smooth backfitting estimator of the $m_j$'s converges at the pointwise rate $n^{-2/5}$ when these functions are twice continuously differentiable. The estimator is asymptotically normally distributed and avoids the curse of dimensionality, but extending it to models in which $F$ is not the identity function appears to be quite difficult. Horowitz, Klemelä and Mammen (2002) (hereinafter HKM) discuss optimality properties of a variety of estimators for nonparametric additive models without link functions.

Tjøstheim and Auestad (1994), Linton and Nielsen (1995), Chen, Härdle, Linton and Severance-Lossin (1996) and Fan, Härdle and Mammen (1998) have investigated the properties of marginal integration estimators for the case in which $F$ is the identity function. These estimators are based on the observation that when $F$ is the identity function, then $m_1(x^1)$, say, is given up to an additive constant by

$$(1.2) \qquad \int \mathbf{E}(Y|X=x)w(x^2,\ldots,x^d)\,dx^2\cdots dx^d,$$

where $w$ is a nonnegative function satisfying

$$\int w(x^2,\ldots,x^d)\,dx^2\cdots dx^d = 1.$$

Therefore, $m_1(x^1)$ can be estimated up to an additive constant by replacing $\mathbf{E}(Y|X=x)$ in (1.2) with a nonparametric estimator. Linton and Nielsen (1995), Chen, Härdle, Linton and Severance-Lossin (1996) and Fan, Härdle and Mammen (1998) have given conditions under which a variety of estimators based on the marginal integration idea converge at rate $n^{-2/5}$ and are asymptotically normal. The latter two estimators have the oracle property. That is, the asymptotic distribution of the estimator of each additive component is the same as it would be if the other components were known. LH extend marginal integration to the case in which $F$ is not the identity function. However, marginal integration estimators have a curse of dimensionality: the smoothness of the $m_j$'s must increase as the dimension of $X$ increases to achieve $n^{-2/5}$ convergence. The reason for this is that estimating $\mathbf{E}(Y|X=x)$ requires carrying out a $d$-dimensional nonparametric regression. If $d$ is large and the $m_j$'s are only twice differentiable, then the bias of the resulting estimator of $\mathbf{E}(Y|X=x)$ converges to zero too slowly as $n \to \infty$ to estimate the $m_j$'s with an $n^{-2/5}$ rate. For example, the estimator of Fan, Härdle and Mammen (1998), which imposes the weakest smoothness conditions of any existing marginal integration estimator, requires more than two derivatives if $d \geq 5$.

This paper describes a two-stage estimation procedure that does not require a $d$-dimensional nonparametric regression and, thereby, avoids the curse of dimensionality. In the first stage, nonlinear least squares is used to obtain a series approximation to each $m_j$. The first-stage procedure imposes the additive structure of (1.1) and yields estimates of the $m_j$'s that have smaller asymptotic biases than do estimators based on marginal integration or other procedures that require $d$-dimensional nonparametric estimation. The first-stage estimates are inputs to the second stage. The second-stage estimate of, say, $m_1$ is obtained by taking one Newton step from the first-stage estimate toward a local linear estimate. In large samples, the second-stage estimator has a structure similar to that of a local linear estimator, so deriving its pointwise rate of convergence and asymptotic distribution is relatively easy. The main results of this paper can also be obtained by using a local constant estimate in the second stage, and the results of Monte Carlo experiments described in Section 5 show that a local constant estimator has better finite-sample performance under some conditions. However, a local linear estimator has better boundary behavior and better ability to adapt to nonuniform designs, among other desirable properties [Fan and Gijbels (1996)].

Our approach differs from typical two-stage estimation, which aims at estimating one unknown parameter or function [e.g., Fan and Chen (1999)]. In this setting, a consistent estimator is obtained in the first stage and is updated in the second, possibly by taking a Newton step toward the optimum of an appropriate objective function. In contrast, in our setting, there are several unknown functions but we update the estimator of only one. It is essential that the first-stage estimators of the other functions have negligible bias. The variances of these estimators must also converge to zero but can have relatively slow rates. We show that asymptotically, the estimation error of the other functions does not appear in the updated estimator of the function of interest.

HKM use a two-stage estimation approach that is similar to the one used here, but HKM do not consider models with link functions, and they use backfitting for the first-stage estimator. Derivation of the properties of a backfitting estimator for a model with a link function appears to be very complicated. We conjecture that a classical backfitting estimator would have the same asymptotic variance as the one in this paper but a different and, possibly, complicated bias. We also conjecture that a classical backfitting estimator would not have the oracle property. Nonetheless, we do not argue here that our procedure outperforms classical backfitting, in the sense of minimizing an optimality criterion such as the asymptotic mean-square error. However, our procedure has the advantages of a complete asymptotic distribution theory and the oracle property.

The remainder of this paper is organized as follows. Section 2 provides an informal description of the two-stage estimator. The main results are presented in Section 3. Section 4 discusses the selection of bandwidths. Section 5 presents the results of a small simulation study, and Section 6 presents concluding comments. The proofs of theorems are in Section 7. Throughout the paper, subscripts index observations and superscripts denote components of vectors. Thus, $X_i$ is the $i$th observation of $X$, $X^j$ is the $j$th component of $X$, and $X_i^j$ is the $i$th observation of the $j$th component.

**2. Informal description of the estimator.** Assume that the support of $X$ is $\mathcal{X} \equiv [-1, 1]^d$, and normalize $m_1, \ldots, m_d$ so that

$$\int_{-1}^{1} m_j(v)\, dv = 0, \qquad j = 1, \ldots, d.$$

For any $x \in \mathbb{R}^d$ define $m(x) = m_1(x^1) + \cdots + m_d(x^d)$, where $x^j$ is the $j$th component of $x$. Let $\{p_k : k = 1, 2, \ldots\}$ denote a basis for smooth functions on $[-1, 1]$. A precise definition of "smooth" and conditions that the basis functions must satisfy are given in Section 3. These conditions include

$$(2.1) \qquad \int_{-1}^{1} p_k(v)\, dv = 0,$$

$$(2.2) \qquad \int_{-1}^{1} p_j(v) p_k(v)\, dv = \begin{cases} 1, & \text{if } j = k, \\ 0, & \text{otherwise,} \end{cases}$$

and

$$(2.3) \qquad m_j(x^j) = \sum_{k=1}^{\infty} \theta_{jk} p_k(x^j),$$

for each $j = 1, \ldots, d$, each $x^j \in [0, 1]$ and suitable coefficients $\{\theta_{jk}\}$. For any positive integer $\kappa$, define

$$P_\kappa(x) = [1, p_1(x^1), \ldots, p_\kappa(x^1), p_1(x^2), \ldots, p_\kappa(x^2), \ldots, p_1(x^d), \ldots, p_\kappa(x^d)]'.$$

Then for $\theta_\kappa \in \mathbb{R}^{\kappa d + 1}$, $P_\kappa(x)'\theta_\kappa$ is a series approximation to $\mu + m(x)$. Section 3 gives conditions that $\kappa$ must satisfy. These require that $\kappa \to \infty$ at an appropriate rate as $n \to \infty$.

To obtain the first-stage estimators of the $m_j$'s, let $\{Y_i, X_i : i = 1, \ldots, n\}$ be a random sample of $(Y, X)$. Let $\hat{\theta}_{n\kappa}$ be a solution to

$$\underset{\theta \in \Theta_\kappa}{\text{minimize}}: \quad S_{n\kappa}(\theta) \equiv n^{-1} \sum_{i=1}^{n} \{Y_i - F[P_\kappa(X_i)'\theta]\}^2,$$

where $\Theta_\kappa \subset \mathbb{R}^{\kappa d + 1}$ is a compact parameter set. The series estimator of $\mu + m(x)$ is

$$\tilde{\mu} + \tilde{m}(x) = P_\kappa(x)'\hat{\theta}_{n\kappa},$$

where $\tilde{\mu}$ is the first component of $\hat{\theta}_{n\kappa}$. The estimator of $m_j(x^j)$ for any $j = 1, \ldots, d$ and any $x^j \in [0, 1]$ is the product of $[p_1(x^j), \ldots, p_\kappa(x^j)]$ with the appropriate components of $\hat{\theta}_\kappa$.

To obtain the second-stage estimator of (say) $m_1(x^1)$, let $\tilde{X}_i$ denote the $i$th observation of $\tilde{X} \equiv (X^2, \ldots, X^d)$. Define $\tilde{m}_{-1}(\tilde{X}_i) = \tilde{m}_2(X_i^2) + \cdots + \tilde{m}_d(X_i^d)$, where $X_i^j$ is the $i$th observation of the $j$th component of $X$ and $\tilde{m}_j$ is the series estimator of $m_j$. Let $K$ be a probability density function on $[-1, 1]$, and define $K_h(v) = K(v/h)$ for any real, positive constant $h$. Conditions that $K$ and $h$ must satisfy are given in Section 3. These include $h \to 0$ at an appropriate rate as $n \to \infty$. Define

$$S'_{nj1}(x^1, \tilde{m}) = -2 \sum_{i=1}^{n} \{Y_i - F[\tilde{\mu} + \tilde{m}_1(x^1) + \tilde{m}_{-1}(\tilde{X}_i)]\}$$
$$\times F'[\tilde{\mu} + \tilde{m}_1(x^1) + \tilde{m}_{-1}(\tilde{X}_i)](X_i^1 - x^1)^j K_h(x^1 - X_i^1)$$

for $j = 0, 1$ and

$$S''_{nj1}(x^1, \tilde{m}) = 2 \sum_{i=1}^{n} F'[\tilde{\mu} + \tilde{m}_1(x^1) + \tilde{m}_{-1}(\tilde{X}_i)]^2 (X_i^1 - x^1)^j K_h(x^1 - X_i^1)$$
$$- 2 \sum_{i=1}^{n} \{Y_i - F[\tilde{\mu} + \tilde{m}_1(x^1) + \tilde{m}_{-1}(\tilde{X}_i)]\}$$
$$\times F''[\tilde{\mu} + \tilde{m}_1(x^1) + \tilde{m}_{-1}(\tilde{X}_i)](X_i^1 - x^1)^j K_h(x^1 - X_i^1)$$

for $j = 0, 1, 2$. The second-stage estimator of $m_1(x^1)$ is

$$\hat{m}_1(x^1) = \tilde{m}_1(x^1) - \frac{S''_{n21}(x^1, \tilde{m}) S'_{n01}(x^1, \tilde{m}) - S''_{n11}(x^1, \tilde{m}) S'_{n11}(x^1, \tilde{m})}{S''_{n01}(x^1, \tilde{m}) S''_{n21}(x^1, \tilde{m}) - S''_{n11}(x^1, \tilde{m})^2}.$$

(2.4)

The second-stage estimators of $m_2(x^2), \ldots, m_d(x^d)$ are obtained similarly. Section 3.3 describes a weighted version of this estimator that minimizes the asymptotic variance of $n^{2/5}[\hat{m}_1(x^1) - m(x^1)]$. However, due to interactions between the weight function and the bias, the weighted estimator does not necessarily minimize the asymptotic mean-square error.

The estimator (2.4) can be understood intuitively as follows. If $\tilde{\mu}$ and $\tilde{m}_{-1}$ were the true values of $\mu$ and $m_{-1}$, the local linear estimator of $m_1(x^1)$ would minimize

$$S_{n1}(x^1, b_0, b_1) = \sum_{i=1}^{n} \{Y_i - F[\tilde{\mu} + b_0 + b_1(X_i^1 - x^1)$$

(2.5)

$$+ \tilde{m}_{-1}(\tilde{X}_i)]\}^2 K_h(x^1 - X_i^1).$$

Moreover, $S'_{nj1}(x^1, \tilde{m}) = \partial S_{n1}(x^1, b_0, b_1)/\partial b_j$ $(j = 0, 1)$ evaluated at $b_0 = \tilde{m}_1(x^1)$ and $b_1 = 0$. $S''_{nj1}(x^1, \tilde{m})$ gives the second derivatives of $S_{n1}(x^1, b_0, b_1)$ evaluated at the same point. The estimator (2.4) is the result of taking one Newton step from the starting values $b_0 = \tilde{m}_1(x^1)$, $b_1 = 0$ toward the minimum of the right-hand side of (2.5).

Section 3 gives conditions under which $\hat{m}_1(x^1) - m_1(x^1) = O_p(n^{-2/5})$ and $n^{2/5}[\hat{m}_1(x^1) - m_1(x^1)]$ is asymptotically normally distributed for any finite $d$ when $F$ and the $m_j$'s are twice continuously differentiable.

**3. Main results.** This section has three parts. Section 3.1 states the assumptions that are used to prove the main results. Section 3.2 states the results. The main results are the $n^{-2/5}$-consistency and asymptotic normality of the $m_j$'s. Section 3.3 describes the weighted estimator.

The following additional notation is used. For any matrix $A$, define the norm $\|A\| = [\text{trace}(A'A)]^{1/2}$. Define $U = Y - F[\mu + m(X)]$, $V(x) = \text{Var}(U|X = x)$, $Q_\kappa = \mathbf{E}\{F'[\mu + m(X)]^2 P_\kappa(X) P_\kappa(X)'\}$, and $\Psi_\kappa = Q_\kappa^{-1} \mathbf{E}\{F'[\mu + m(X)]^2 V(X) \times P_\kappa(X) P_\kappa(X)'\} Q_\kappa^{-1}$ whenever the latter quantity exists. $Q_\kappa$ and $\Psi_\kappa$ are $d(\kappa) \times d(\kappa)$ positive semidefinite matrices, where $d(\kappa) = \kappa d + 1$. Let $\lambda_{\kappa,\min}$ denote the smallest eigenvalue of $Q_\kappa$. Let $Q_{\kappa,ij}$ denote the $(i, j)$ element of $Q_\kappa$. Define $\zeta_\kappa = \sup_{x \in \mathcal{X}} \|P_\kappa(x)\|$. Let $\{\theta_{jk}\}$ be the coefficients of the series expansion (2.3). For each $\kappa$ define

$$\theta_\kappa = (\mu, \theta_{11}, \ldots, \theta_{1\kappa}, \theta_{21}, \ldots, \theta_{2\kappa}, \ldots, \theta_{d1}, \ldots, \theta_{d\kappa})'.$$

3.1. *Assumptions.* The main results are obtained under the following assumptions.

ASSUMPTION A1. The data, $\{(Y_i, X_i) : i = 1, \ldots, n\}$, are an i.i.d. random sample from the distribution of $(Y, X)$, and $\mathbf{E}(Y|X = x) = F[\mu + m(x)]$ for almost every $x \in \mathcal{X} \equiv [-1, 1]^d$.

ASSUMPTION A2. (i) The support of $X$ is $\mathcal{X}$.

(ii) The distribution of $X$ is absolutely continuous with respect to Lebesgue measure.

(iii) The probability density function of $X$ is bounded, bounded away from zero and twice continuously differentiable on $\mathcal{X}$.

(iv) There are constants $c_V > 0$ and $C_V < \infty$ such that $c_V \leq \mathrm{Var}(U|X = x) \leq C_V$ for all $x \in \mathcal{X}$.

(v) There is a constant $C_U < \infty$ such that $\mathbf{E}|U|^j \leq C_U^{j-2} j! \mathbf{E}(U^2) < \infty$ for all $j \geq 2$.

ASSUMPTION A3. (i) There is a constant $C_m < \infty$ such that $|m_j(v)| \leq C_m$ for each $j = 1, \ldots, d$ and all $v \in [-1, 1]$.

(ii) Each function $m_j$ is twice continuously differentiable on $[-1, 1]$.

(iii) There are constants $C_{F1} < \infty, c_{F2} > 0$, and $C_{F2} < \infty$ such that $F(v) \leq C_{F1}$ and $c_{F2} \leq F'(v) \leq C_{F2}$ for all $v \in [\mu - C_m d, \mu + C_m d]$.

(iv) $F$ is twice continuously differentiable on $[\mu - C_m d, \mu + C_m d]$.

(v) There is a constant $C_{F3} < \infty$ such that $|F''(v_2) - F''(v_1)| \leq C_{F3}|v_2 - v_1|$ for all $v_2, v_1 \in [\mu - C_m d, \mu + C_m d]$.

ASSUMPTION A4. (i) There are constants $C_Q < \infty$ and $c_\lambda > 0$ such that $|Q_{\kappa, ij}| \leq C_Q$ and $\lambda_{\kappa,\min} > c_\lambda$ for all $\kappa$ and all $i, j = 1, \ldots, d(\kappa)$.

(ii) The largest eigenvalue of $\Psi_\kappa$ is bounded for all $\kappa$.

ASSUMPTION A5. (i) The functions $\{p_k\}$ satisfy (2.1) and (2.2).

(ii) There is a constant $c_\kappa > 0$ such that $\zeta_\kappa \geq c_\kappa$ for all sufficiently large $\kappa$.

(iii) $\zeta_\kappa = O(\kappa^{1/2})$ as $\kappa \to \infty$.

(iv) There are a constant $C_\theta < \infty$ and vectors $\theta_{\kappa 0} \in \Theta_\kappa \equiv [-C_\theta, C_\theta]^{d(\kappa)}$ such that $\sup_{x \in \mathcal{X}} |\mu + m(x) - P_\kappa(x)'\theta_{\kappa 0}| = O(\kappa^{-2})$ as $\kappa \to \infty$.

(v) For each $\kappa, \theta_\kappa$ is an interior point of $\Theta_\kappa$.

ASSUMPTION A6. (i) $\kappa = C_\kappa n^{4/15+\nu}$ for some constant $C_\kappa$ satisfying $0 < C_\kappa < \infty$ and some $\nu$ satisfying $0 < \nu < 1/30$.

(ii) $h = C_h n^{-1/5}$ for some constant $C_h$ satisfying $0 < C_h < \infty$.

ASSUMPTION A7. The function $K$ is a bounded, continuous probability density function on $[-1, 1]$ and is symmetric about $0$.

The assumption that the support of $X$ is $[-1, 1]^d$ entails no loss of generality as it can always be satisfied by carrying out monotone increasing transformations of the components of $X$, even if their support before transformation is unbounded. For practical computations, it suffices to transform the empirical support to $[-1, 1]^d$. Assumption A2 precludes the possibility of treating discrete covariates with our method, though they can be handled inelegantly by conditioning on them. Another possibility is to develop a version of our estimator for a partially linear generalized additive model in which discrete covariates are included in the parametric (linear) term. However, this extension is beyond the scope of the present paper. Differentiability of the density of $X$ [Assumption A2(iii)] is used to insure that the bias of our estimator converges to zero sufficiently rapidly. Assumption A2(v) restricts the thickness of the tails of the distribution of $U$ and is used to prove consistency of the first-stage estimator. Assumption A3 defines the sense in which $F$ and the $m_j$'s must be smooth. Assumption A3(iii) is needed for identification. Assumption A4 insures the existence and nonsingularity of the covariance matrix of the asymptotic form of the first-stage estimator. This is analogous to assuming that the information matrix is positive definite in parametric maximum likelihood estimation. Assumption A4(i) implies Assumption A4(ii) if $U$ is homoskedastic. Assumptions A5(iii) and A5(iv) bound the magnitudes of the basis functions and insure that the errors in the series approximations to the $m_j$'s converge to zero sufficiently rapidly as $\kappa \to \infty$. These assumptions are satisfied by spline and (for periodic functions) Fourier bases. Assumption A6 states the rates at which $\kappa \to \infty$ and $h \to 0$ as $n \to \infty$. The assumed rate of convergence of $h$ is well known to be asymptotically optimal for one-dimensional kernel mean-regression when the conditional mean function is twice continuously differentiable. The required rate for $\kappa$ insures that the asymptotic bias and variance of the first-stage estimator are sufficiently small to achieve an $n^{-2/5}$ rate of convergence in the second stage. The $L_2$ rate of convergence of a series estimator of $m_j$ is maximized by setting $\kappa \propto n^{1/5}$, which is slower than the rates permitted by Assumption A6(i) [Newey (1997)]. Thus, Assumption A6(i) requires the first-stage estimator to be undersmoothed. Undersmoothing is needed to insure sufficiently rapid convergence of the bias of the first-stage estimator. We show that the first-order performance of our second-stage estimator does not depend on the choice of $\kappa$ if Assumption A6(i) is satisfied. See Theorems 2 and 3. Optimizing the choice of $\kappa$ would require a rather complicated higher-order theory and is beyond the scope of this paper, which is restricted to first-order asymptotics.

3.2. *Theorems.* This section states two theorems that give the main results of the paper. Theorem 1 gives the asymptotic behavior of the first-stage series estimator under Assumptions A1–A6(i). Theorem 2 gives the properties of the second-stage estimator. For $i = 1, \ldots, n$, define $U_i = Y_i - F[\mu + m(X_i)]$ and $b_{\kappa 0}(x) = \mu + m(x) - P_\kappa(x)'\theta_{\kappa 0}$. Let $\|v\|$ denote the Euclidean norm of any finite-dimensional vector $v$.

THEOREM 1. *Let Assumptions* A1–A6(i) *hold. Then:*

(a) $\lim_{n \to \infty} \|\hat{\theta}_{n\kappa} - \theta_{\kappa 0}\| = 0$ *almost surely,*

(b) $\hat{\theta}_{n\kappa} - \theta_{\kappa 0} = O_p(\kappa^{1/2}/n^{1/2} + \kappa^{-2})$, *and*

(c) $\sup_{x \in \mathcal{X}} |\tilde{m}(x) - m(x)| = O_p(\kappa/n^{1/2} + \kappa^{-3/2})$.

*In addition:*

(d) $\hat{\theta}_{n\kappa} - \theta_{\kappa 0} = n^{-1} Q_\kappa^{-1} \sum_{i=1}^n F'[\mu + m(X_i)] P_\kappa(X_i) U_i + n^{-1} Q_\kappa^{-1} \times \sum_{i=1}^n F'[\mu + m(X_i)]^2 P_\kappa(X_i) b_{\kappa 0}(X_i) +$ $R_n$, *where* $\|R_n\| = O_p(\kappa^{3/2}/n + n^{-1/2})$.

Now let $f_X$ denote the probability density function of $X$. For $j = 0, 1$, define

$$S'_{nj1}(x^1, m) = -2 \sum_{i=1}^n \{Y_i - F[\mu + m_1(x^1) + m_{-1}(\tilde{X}_i)]\}$$
$$\times F'[\mu + m_1(x^1) + m_{-1}(\tilde{X}_i)](X_i^1 - x^1)^j K_h(x^1 - X_i^1).$$

Also define

$$D_0(x^1) = 2 \int F'[\mu + m_1(x^1) + m_{-1}(\tilde{x})]^2 f_X(x^1, \tilde{x}) \, d\tilde{x},$$

$$D_1(x^1) = 2 \int F'[\mu + m_1(x^1) + m_{-1}(\tilde{x})]^2 [\partial f_X(x^1, \tilde{x})/\partial x^1] \, d\tilde{x},$$

$$A_K = \int_{-1}^1 v^2 K(v) \, dv,$$

$$B_K = \int_{-1}^1 K(v)^2 \, dv,$$

$$g(x^1, \tilde{x}) = F''[\mu + m_1(x^1) + m_{-1}(\tilde{x})] m_1'(x^1)$$
$$+ F'[\mu + m_1(x^1) + m_{-1}(\tilde{x})] m_1''(x^1),$$

$$\beta_1(x^1) = 2 C_h^2 A_K D_0(x^1)^{-1}$$
$$\times \int g(x^1, \tilde{x}) F'[\mu + m_1(x^1) + m_{-1}(\tilde{x})] f_X(x^1, \tilde{x}) \, d\tilde{x}$$

and

$$V_1(x^1) = B_K C_h^{-1} D_0(x^1)^{-2}$$
$$\times \int \text{Var}(U|x^1, \tilde{x}) F'[\mu + m_1(x^1) + m_{-1}(\tilde{x})]^2 f_X(x^1, \tilde{x}) \, d\tilde{x}.$$

The next theorem gives the asymptotic properties of the second-stage estimator.

THEOREM 2. *Let Assumptions* A1–A6 *hold. Then:*

(a) $\hat{m}_1(x^1) - m_1(x^1) = [nhD_0(x^1)]^{-1}\{-S'_{n01}(x^1, m) + [D_1(x^1)/D_0(x^1)] \times S'_{n11}(x^1, m)\} + o_p(n^{-2/5})$ *uniformly over* $|x^1| \le 1 - h$ *and* $\hat{m}_1(x^1) - m_1(x^1) = O_p[(\log n)^{1/2}n^{-2/5}]$ *uniformly over* $|x^1| \le 1$.

(b) $n^{2/5}[\hat{m}_1(x^1) - m_1(x^1)] \xrightarrow{d} N[\beta_1(x^1), V_1(x^1)]$.

(c) *If* $j \ne 1$, *then* $n^{2/5}[\hat{m}_1(x^1) - m_1(x^1)]$ *and* $n^{2/5}[\hat{m}_j(x^j) - m_j(x^j)]$ *are asymptotically independently normally distributed.*

Theorem 2(a) implies that asymptotically, $n^{2/5}[\hat{m}_1(x^1) - m_1(x^1)]$ is not affected by random sampling errors in the first-stage estimator. In fact, the second-stage estimator of $m_1(x^1)$ has the same asymptotic distribution that it would have if $m_2, \ldots, m_d$ were known and local-linear estimation were used to estimate $m_1(x^1)$ directly. In this sense, our estimator has an oracle property. Parts (b) and (c) of Theorem 2 imply that the estimators of $m_1(x^1), \ldots, m_d(x^d)$ are asymptotically independently distributed.

It is also possible to use a local-constant estimator in the second stage. The resulting second-stage estimator is

$$\hat{m}_{1,LC}(x^1) = \tilde{m}_1(x^1) - S'_{n01}(x^1, \tilde{m})/S''_{n01}(x^1, \tilde{m}).$$

The following modification of Theorem 2, which we state without proof, gives the asymptotic properties of the local-constant second-stage estimator. Define

$$g_{LC}(x^1, \tilde{x}) = (\partial^2/\partial\zeta^2)\{F[m_1(\zeta + x^1) + m_{-1}(\tilde{x})]$$
$$- F[m_1(x^1) + m_{-1}(\tilde{x})]\}f_X(\zeta + x^1, \tilde{x})|_{\zeta=0}$$

and

$$\beta_{1,LC}(x^1) = 2C_h^2 A_K D_0(x^1)^{-1}$$
$$\times \int g_{LC}(x^1, \tilde{x})F'[\mu + m_1(x^1) + m_{-1}(\tilde{x})]f_X(x^1, \tilde{x})\, d\tilde{x}.$$

THEOREM 3. *Let Assumptions* A1–A6 *hold. Then*

(a) $\hat{m}_{1,LC}(x^1) - m_1(x^1) = -[nhD_0(x^1)]^{-1}S'_{n01}(x^1, m) + o_p(n^{-2/5})$ *uniformly over* $|x^1| \le 1 - h$ *and* $\hat{m}_1(x^1) - m_1(x^1) = O_p[(\log n)^{1/2}n^{-2/5}]$ *uniformly over* $|x^1| \le 1$.

(b) $n^{2/5}[\hat{m}_{1,LC}(x^1) - m_1(x^1)] \xrightarrow{d} N[\beta_{1,LC}(x^1), V_1(x^1)]$.

(c) *If $j \neq 1$, then $n^{2/5}[\hat{m}_{1,LC}(x^1) - m_1(x^1)]$ and $n^{2/5}[\hat{m}_{j,LC}(x^j) - m_j(x^j)]$ are asymptotically independently normally distributed.*

$V_1(x^1)$ and $\beta_1(x^1)$ and $\beta_{1,LC}(x^1)$ can be estimated consistently by replacing unknown population parameters with consistent estimators. Section 4 gives a method for estimating the derivatives of $m_1$ that are in the expressions for $\beta_1(x^1)$ and $\beta_{1,LC}(x^1)$. As is usual in nonparametric estimation, reasonably precise bias estimation is possible only by making assumptions that amount to undersmoothing. One way of doing this is to assume that the second derivative of $m_1$ satisfies a Lipschitz condition. Alternatively, one can set $h = C_h n^{-\gamma}$ for $1/5 < \gamma < 1$. Then $n^{(1-\gamma)/2}[\hat{m}_1(x^1) - m_1(x^1)] \xrightarrow{d} N[0, V_1(x^1)]$, and $n^{(1-\gamma)/2}[\hat{m}_{1,LC}(x^1) - m_1(x^1)] \xrightarrow{d} N[0, V_1(x^1)]$.

**3.3. *A weighted estimator.*** A weighted estimator can be obtained by replacing $S'_{nj1}(x^1, \tilde{m})$ and $S''_{nj1}(x^1, \tilde{m})$ in (2.5) with

$$S'_{nj1}(x^1, \tilde{m}, w) = -2 \sum_{i=1}^{n} w(x^1, \tilde{X}_i)\{Y_i - F[\tilde{\mu} + \tilde{m}_1(x^1) + \tilde{m}_{-1}(\tilde{X}_i)]\}$$
$$\times F'[\tilde{\mu} + \tilde{m}_1(x^1) + \tilde{m}_{-1}(\tilde{X}_i)](X_i^1 - x^1)^j K_h(x^1 - X_i^1)$$

and

$$S''_{nj1}(x^1, \tilde{m}, w)$$
$$= 2 \sum_{i=1}^{n} w(x^1, \tilde{X}_i) F'[\tilde{\mu} + \tilde{m}_1(x^1) + \tilde{m}_{-1}(\tilde{X}_i)]^2 (X_i^1 - x^1)^j K_h(x^1 - X_i^1)$$
$$- 2 \sum_{i=1}^{n} w(x^1, \tilde{X}_i)\{Y_i - F[\tilde{\mu} + \tilde{m}_1(x^1) + \tilde{m}_{-1}(\tilde{X}_i)]\}$$
$$\times F''[\tilde{\mu} + \tilde{m}_1(x^1) + \tilde{m}_{-1}(\tilde{X}_i)](X_i^1 - x^1)^j K_h(x^1 - X_i^1)$$

for $j = 0, 1, 2$, where $w$ is a nonnegative weight function that is assumed for the moment to be nonstochastic. It is convenient to normalize $w$ so that

$$\int w(x^1, \tilde{x}) F'[\mu + m_1(x^1) + m_{-1}(\tilde{x})]^2 f_X(x^1, \tilde{x}) \, d\tilde{x} = 1$$

for each $x^1 \in [-1, 1]$. Arguments identical to those used to prove Theorem 2 show that the variance of the asymptotic distribution of the resulting local-linear or local-constant estimator of $m_1(x^1)$ is

$$V_1(x^1, w) = 0.25 B_K C_h^{-1} \int w(x^1, \tilde{x})^2 \operatorname{Var}(U|x^1, \tilde{x})$$
$$\times F'[\mu + m_1(x^1) + m_{-1}(\tilde{x})]^2 f_X(x^1, \tilde{x}) \, d\tilde{x}.$$

It follows from Lemma 1 of Fan, Härdle and Mammen (1998) that $V(x^1, w)$ is minimized by setting $w(x^1, \tilde{x}) \propto 1/\operatorname{Var}(U|x^1, \tilde{x})$, thereby yielding

$$V_1(x^1, w) = 0.25 B_K C_h^{-1} D_2(x^1)^{-1} \int F'[\mu + m_1(x^1) + m_{-1}(\tilde{x})]^2 f_X(x^1, \tilde{x}) \, d\tilde{x},$$

where

$$D_2(x^1) = \int \operatorname{Var}(U|x^1, \tilde{x})^{-1} F'[\mu + m_1(x^1) + m_{-1}(\tilde{x})]^2 f_X(x^1, \tilde{x}) \, d\tilde{x}.$$

In an application, it suffices to replace the variance-minimizing weight function with a consistent estimator. For example, $F'[\mu + m_1(x^1) + m_{-1}(\tilde{x})]$ can be estimated from the first estimation stage, $\operatorname{Var}(U|x^1, \tilde{x})$ can be estimated by applying a nonparametric regression to the squared residuals of the first-stage estimate and kernel methods can be used to estimate $f_X(x^1, \tilde{x})$.

The minimum-variance estimator is not a minimum asymptotic mean-square error estimator unless undersmoothing is used to remove the asymptotic bias of $\hat{m}_1$. This is because weighting affects the bias when the latter is nonnegligible. The weight function that minimizes the asymptotic mean-square error is the solution to an integral equation and does not have a closed-form analytic representation.

**4. Bandwidth selection.** This section presents a plug-in and a penalized least squares (PLS) method for choosing $h$ in applications. We begin with a description of the plug-in method. This method estimates the value of $h$ that minimizes the asymptotic integrated mean-square error (AIMSE) of $n^{2/5}[\hat{m}_1(x^1) - m_1(x^1)]$ for $j = 1, \ldots, d$. We discuss only local-linear estimation, but similar results hold for local-constant estimation. The AIMSE of $n^{2/5}(\hat{m}_1 - m_1)$ is defined as

$$\text{AIMSE}_1 = n^{4/5} \int_{-1}^{1} w(x^1)[\beta_1(x^1)^2 + V_1(x^1)] \, dx^1,$$

where $w(\cdot)$ is a nonnegative weight function that integrates to 1. We also define the integrated squared error (ISE) as

$$\text{ISE}_1 = n^{4/5} \int_{-1}^{1} w(x^1)[\hat{m}_1(x^1) - m_1(x^1)]^2 \, dx^1.$$

We define the asymptotically optimal bandwidth for estimating $m_1$ as $C_{h1} n^{-1/5}$, where $C_{h1}$ minimizes $\text{AIMSE}_1$. Let

$$\tilde{\beta}_1(x^1) = \beta_1(x^1)/C_h^2 \quad \text{and} \quad \tilde{V}_1(x^1) = C_h V_1(x^1).$$

Then

$$(4.1) \qquad C_{h1} = \left[ \frac{1}{4} \frac{\int_{-1}^{1} w(x^1) \tilde{V}_1(x^1) \, dx^1}{\int_{-1}^{1} w(x^1) \tilde{\beta}_1(x^1)^2 \, dx^1} \right]^{1/5}.$$

The results for the plug-in method rely on the following two theorems. Theorem 4 shows that the difference between the ISE and AIMSE is asymptotically negligible. Theorem 5 gives a method for estimating the first and second derivatives of $m_j$. Let $G^{(\ell)}$ denote the $\ell$th derivative of any $\ell$-times differentiable function $G$.

THEOREM 4. *Let Assumptions* A1–A6 *hold. Then for a continuous weight function* $w(\cdot)$ *and as* $n \to \infty$, $\mathrm{AIMSE}_1 = \mathrm{ISE}_1 + o_p(1)$.

THEOREM 5. *Let Assumptions* A1–A6 *hold. Let $L$ be a twice differentiable probability density function on* $[-1, 1]$, *and let* $\{g_n : n = 1, 2, \dots\}$ *be a sequence of strictly positive real numbers satisfying* $g_n \to 0$ *and* $g_n^2 n^{4/5} (\log n)^{-1} \to \infty$ *as* $n \to \infty$. *For* $\ell = 1, 2$ *define*

$$\hat{m}_1^{(\ell)}(x^1) = g_n^{-1-\ell} \int_{-1}^{1} L^{(\ell)}[(x^1 - v)/g_n] \hat{m}_1(v) \, dv.$$

*Then as* $n \to \infty$ *and for* $\ell = 1, 2$,

$$\sup_{|x^1| \leq 1} |\hat{m}_1^{(\ell)}(x^1) - m_1^{(\ell)}(x^1)| = o_p(1).$$

A plug-in estimator of $C_{h1}$ can now be obtained by replacing unknown population quantities on the right-hand side of (4.1) with consistent estimators. Theorem 5 provides consistent estimators of the required derivatives of $m_1$. Estimators of the conditional variance of $U$ and of $f_X$ can be obtained by using standard kernel methods.

We now describe the PLS method. This method simultaneously estimates the bandwidths for second-stage estimation of all of the functions $m_j$ ($j = 1, \dots, d$). Let $h_j = C_{hj} n^{-1/5}$ be the bandwidth for $\hat{m}_j$. Then the PLS method selects the $C_{hj}$'s that minimize an estimate of the average squared error (ASE),

$$\mathrm{ASE}(\bar{h}) = n^{-1} \sum_{i=1}^{n} \{F[\tilde{\mu} + \hat{m}(X_i)] - F[\mu + m(X_i)]\}^2,$$

where $\bar{h} = (C_{h1} n^{-1/5}, \dots, C_{hd} n^{-1/5})$. Specifically, the PLS method selects the $C_{hj}$'s to

$$
\underset{C_{h1}, \dots, C_{hd}}{\text{minimize:}} \quad \mathrm{PLS}(\bar{h}) = n^{-1} \sum_{i=1}^{n} \{Y_i - F[\tilde{\mu} + \hat{m}(X_i)]\}^2
$$

$$
(4.2) \qquad\qquad\qquad + 2K(0) n^{-1} \sum_{i=1}^{n} \{F'[\mu + \hat{m}(X_i)]^2 \hat{V}(X_i)\}
$$

$$
\times \sum_{j=1}^{d} [n^{4/5} C_{hj} \hat{D}_j(X_i^j)]^{-1},
$$

where the $C_{hj}$'s are restricted to a compact, positive interval that excludes 0,

$$\hat{D}_j(x^j) = \frac{1}{nh_j} \sum_{i=1}^{n} K_{h_j}(X_i^j - x^j) F'[\tilde{\mu} + \hat{m}(X_i)]^2$$

and

$$\hat{V}(x) = \left[ \sum_{i=1}^{n} K_{h_1}(X_i^1 - x^1) \cdots K_{h_d}(X_i^d - x^d) \right]^{-1}$$

$$\times \sum_{i=1}^{n} K_{h_1}(X_i^1 - x^1) \cdots K_{h_d}(X_i^d - x^d)\{Y_i - F[\tilde{\mu} + \hat{m}(X_i)]\}^2.$$

The bandwidths used for $\hat{V}$ may be different from those used for $\hat{m}$ because $\hat{V}$ is a full-dimensional nonparametric estimator. We now argue that the difference

$$n^{-1} \sum_{i=1}^{n} U_i^2 + \text{ASE}(\bar{h}) - \text{PLS}(\bar{h})$$

is asymptotically negligible and, therefore, that the solution to (4.2) estimates the bandwidths that minimize ASE. A proof of this result only requires additional smoothness conditions on $F$ and more restrictive assumptions on $\kappa$. The proof can be carried out by making arguments similar to those used in the proof of Theorem 2 but with a higher-order stochastic expansion for $\hat{m} - m$. Here, we provide only a heuristic outline. For this purpose, note that

$$n^{-1} \sum_{i=1}^{n} U_i^2 + \text{ASE}(\bar{h}) - \text{PLS}(\bar{h})$$

$$= 2n^{-1} \sum_{i=1}^{n} \{F[\tilde{\mu} + \hat{m}(X_i)] - F[\mu + m(X_i)]\} U_i$$

$$- 2K(0)n^{-1} \sum_{i=1}^{n} F'[\mu + m(X_i)]^2 \hat{V}(X_i) \sum_{j=1}^{d} [n^{4/5} C_{hj} \hat{D}_j(X_i^j)]^{-1}.$$

We now approximate $F[\tilde{\mu} + \hat{m}(X_i)] - F[\mu + m(X_i)]$ by a linear expansion in $\hat{m} - m$ and replace $\hat{m} - m$ with the stochastic approximation of Theorem 2(a). (A rigorous argument would require a higher-order expansion of $\hat{m} - m$.) Thus, $F[\tilde{\mu} + \hat{m}(X_i)] - F[\mu + m(X_i)]$ is approximated by a linear form in $U_i$. Dropping higher-order terms leads to an approximation of

$$\frac{2}{n} \sum_{i=1}^{n} \{F[\tilde{\mu} + \hat{m}(X_i)] - F[\mu + m(X_i)]\} U_i$$

that is a $U$ statistic in $U_i$. The off-diagonal terms of the $U$ statistic can be shown to be of higher order and, therefore, asymptotically negligible. Thus, we get

$$\frac{2}{n}\sum_{i=1}^{n}\{F[\tilde{\mu}+\hat{m}(X_i)]-F[\mu+m(X_i)]\}U_i$$

$$\approx \frac{2}{n}\sum_{i=1}^{n}F'[\mu+m(X_i)]^2\operatorname{Var}(U_i|X_i)\sum_{j=1}^{d}[n^{4/5}C_{hj}D_{0j}(X_i^j)]h^{-1}K(0),$$

where

$$D_{0j}(x^j)=2\mathbf{E}\{F'[\mu+m(X_i)]^2|X_i^j=x^j\}f_{X^j}(x^j)$$

and $f_{X^j}$ is the probability density function of $X^j$. Now by standard kernel smoothing arguments, $D_{0j}(x^j)\approx\hat{D}_j(x^j)$. In addition, it is clear that $\hat{V}(X_i)\approx V(U_i|X_i)$, which establishes the desired result.

**5. Monte Carlo experiments.** This section presents the results of a small set of Monte Carlo experiments that compare the finite-sample performances of the two-stage estimator, the estimator of LH and the infeasible oracle estimator in which all additive components but one are known. The oracle estimator cannot be used in applications but provides a benchmark against which our feasible estimator can be compared. The infeasible oracle estimator was calculated by solving (2.5).

Experiments were carried out with $d=2$ and $d=5$. The sample size is $n=500$. The experiments with $d=2$ consist of estimating $f_1$ and $f_2$ in the binary logit model

$$\mathbf{P}(Y=1|X=x)=L[f_1(x^1)+f_2(x^2)],$$

where $L$ is the cumulative logistic distribution function

$$L(v)=e^v/[1+e^v],\qquad-\infty<v<\infty.$$

The experiments with $d=5$ consist of estimating $f_1$ and $f_2$ in the binary logit model

$$\mathbf{P}(Y=1|X=x)=L\left[f_1(x^1)+f_2(x^2)+\sum_{j=3}^{5}x^j\right].$$

In all of the experiments,

$$f_1(x)=\sin(\pi x)\quad\text{and}\quad f_2(x)=\Phi(3x),$$

where $\Phi$ is the standard normal distribution function. The components of $X$ are independently distributed as $U[-1,1]$. Estimation is carried out under the assumption that the additive components have two (but not necessarily more) continuous derivatives. Under this assumption, the two-stage estimator has the rate of convergence $n^{-2/5}$. The LH estimator has this rate of convergence if $d = 2$ but not if $d = 5$.

$B$-splines were used for the first stage of the two-stage estimator. The kernel used for the second stage and for the LH estimator is

$$K(v) = \tfrac{15}{16}(1 - v^2)^2 I(|v| \leq 1).$$

Experiments were carried out using both local-constant and local-linear estimators in the second stage of the two-stage method. There were 1000 Monte Carlo replications per experiment with the two-stage estimator but only 500 replications with the LH estimator because of the very long computing times it entails. The experiments were carried out in GAUSS using GAUSS random number generators.

The results of the experiments are summarized in Table 1, which shows the empirical integrated mean-square errors (EIMSEs) of the estimators at the values of the tuning parameters that minimize the EIMSEs. Lengthy computing times precluded using data-based methods for selecting tuning parameters in the experiments. The EIMSEs of the local-constant and local-linear two-stage estimates of $f_1$ are considerably smaller than the EIMSEs of the LH estimator. The EIMSEs of the local-constant and LH estimators of $f_2$ are approximately equal whereas the local-linear estimator of $f_2$ has a larger EIMSE. There is little difference between the EIMSEs of the two-stage local-linear and infeasible oracle estimators. This result is consistent with the oracle property of the two-stage estimator.

**6. Conclusions.** This paper has described an estimator of the additive components of a nonparametric additive model with a known link function. The approach is very general and may be applicable to a wide variety of other models. The estimator is asymptotically normally distributed and has a pointwise rate of convergence in probability of $n^{-2/5}$ when the unknown functions are twice continuously differentiable, regardless of the dimension of the explanatory variable $X$. In contrast, achieving the rate of convergence $n^{-2/5}$ with the only other currently available estimator for this model requires the additive components to have an increasing number of derivatives as the dimension of $X$ increases. In addition, the new estimator has an oracle property: the asymptotic distribution of the estimator of each additive component is the same as it would be if the other components were known.

**7. Proofs of theorems.** Assumptions A1–A7 hold throughout this section.

Table 1
*Results of Monte Carlo experiments*[*]

| Estimator | $\kappa_1$ | $\kappa_2$ | $h_1$ | $h_2$ | Empirical IMSE | |
| | | | | | $f_1$ | $f_2$ |
|---|---|---|---|---|---|---|
| | | | $d = 2$ | | | |
| LH | | | 0.9 | 0.9 | 0.116 | 0.015 |
| Two-stage with local-constant smoothing | 2 | 2 | 0.4 | 0.9 | 0.052 | 0.015 |
| Two-stage with local-linear smoothing | 4 | 2 | 0.5 | 1.4 | 0.052 | 0.023 |
| Infeasible oracle estimator | | | 0.6 | 1.7 | 0.056 | 0.021 |
| | | | $d = 5$ | | | |
| LH | | | 1.0 | 1.0 | 0.145 | 0.019 |
| Two-stage with local-constant smoothing | 2 | 2 | 0.4 | 0.9 | 0.060 | 0.018 |
| Two-stage with local-linear smoothing | 2 | 2 | 0.6 | 1.3 | 0.057 | 0.029 |
| Infeasible oracle estimator | | | 0.6 | 2.0 | 0.057 | 0.023 |

[*]In the two-stage estimator, $\kappa_j$ and $h_j$ $(j = 1, 2)$ are the series length and bandwidth used to estimate $f_j$. In the LH estimator, $h_j$ $(j = 1, 2)$ is the bandwidth used to estimate $f_j$. The values of $\kappa_1$, $\kappa_2$, $h_1$ and $h_2$ minimize the IMSEs of the estimates.

7.1. *Theorem* 1. This section begins with lemmas that are used to prove Theorem 1.

LEMMA 1. *There are constants* $a > 0$ *and* $C < \infty$ *such that*

$$\mathbf{P}\left[\sup_{\theta \in \Theta_\kappa} |S_{n\kappa}(\theta) - \mathbf{E}S_{nk}(\theta)| > \varepsilon\right] \leq C \exp(-na\varepsilon^2)$$

*for any sufficiently small* $\varepsilon > 0$ *and all sufficiently large* $n$.

PROOF. Write

$$S_{n\kappa}(\theta) = n^{-1} \sum_{i=1}^n Y_i^2 - 2S_{n\kappa 1}(\theta) + S_{n\kappa 2}(\theta),$$

where

$$S_{n\kappa 1}(\theta) = n^{-1} \sum_{i=1}^n Y_i F[P_\kappa(X_i)'\theta]$$

and

$$S_{n\kappa 2}(\theta) = n^{-1} \sum_{i=1}^{n} F[P_\kappa(X_i)'\theta]^2.$$

It suffices to prove that

$$\mathbf{P}\bigg[\sup_{\theta \in \Theta_\kappa} |S_{n\kappa j}(\theta) - \mathbf{E}S_{n\kappa j}(\theta)| > \varepsilon\bigg] \le \tilde{C}\exp(-na\varepsilon^2) \qquad (j = 1, 2)$$

for any $\varepsilon > 0$, some $\tilde{C} < \infty$ and all sufficiently large $n$. The proof is given only for $j = 1$. Similar arguments apply when $j = 2$.

Define $\tilde{S}_{n\kappa 1}(\theta) = S_{n\kappa 1}(\theta) - \mathbf{E}S_{n\kappa 1}(\theta)$. Divide $\Theta_\kappa$ into hypercubes of edge-length $\ell$. Let $\Theta_\kappa^{(1)}, \dots, \Theta_\kappa^{(M)}$ denote the $M = (2C_\theta/\ell)^{d(\kappa)}$ cubes thus created. Let $\theta_{\kappa j}$ be the point at the center of $\Theta_\kappa^{(j)}$. The maximum distance between $\theta_{\kappa j}$ and any other point in $\Theta_\kappa^{(j)}$ is $r = d(\kappa)^{1/2}\ell/2$, and $M = \exp\{d(\kappa)[\log(C_\theta/r) + (1/2)\log d(\kappa)]\}$. Now

$$\bigg[\sup_{\theta \in \Theta_\kappa} |\tilde{S}_{n\kappa 1}(\theta)| > \varepsilon\bigg] \subset \bigcup_{j=1}^{M}\bigg[\sup_{\theta \in \Theta_\kappa^{(j)}} |\tilde{S}_{n\kappa 1}(\theta)| > \varepsilon\bigg].$$

Therefore,

$$\mathbf{P}_n \equiv \mathbf{P}\bigg[\sup_{\theta \in \Theta_\kappa} |\tilde{S}_{n\kappa 1}(\theta)| > \varepsilon\bigg] \le \sum_{j=1}^{M}\mathbf{P}\bigg[\sup_{\theta \in \Theta_\kappa^{(j)}} |\tilde{S}_{n\kappa 1}(\theta)| > \varepsilon\bigg].$$

Now for $\theta \in \Theta_\kappa^{(j)}$,

$$|\tilde{S}_{n\kappa 1}(\theta)| \le |\tilde{S}_{n\kappa 1}(\theta_{\kappa j})| + |\tilde{S}_{n\kappa 1}(\theta) - \tilde{S}_{n\kappa 1}(\theta_{\kappa j})|$$

$$\le |\tilde{S}_{n\kappa 1}(\theta_{\kappa j})| + C_{F2}\zeta_\kappa r\bigg(n^{-1}\sum_{i=1}^{n}|Y_i| + C_{F1}\bigg)$$

for all sufficiently large $\kappa$ and, therefore, $n$. Therefore, for all sufficiently large $n$,

$$\mathbf{P}\bigg[\sup_{\theta \in \Theta_\kappa^{(j)}} |\tilde{S}_{n\kappa 1}(\theta)| > \varepsilon\bigg]$$

$$\le \mathbf{P}[|\tilde{S}_{n\kappa 1}(\theta_{\kappa j})| > \varepsilon/2] + \mathbf{P}\bigg[C_{F2}\zeta_\kappa r\bigg(n^{-1}\sum_{i=1}^{n}|Y_i| + C_{F1}\bigg) > \varepsilon/2\bigg].$$

Choose $r = \zeta_\kappa^{-2}$. Then $\varepsilon/2 - C_{F2}\zeta_\kappa r[C_{F1} + \mathbf{E}(|Y|)] > \varepsilon/4$ for all sufficiently large $\kappa$. Moreover,

$$\mathbf{P}\bigg[C_{F2}\zeta_\kappa r\bigg(n^{-1}\sum_{i=1}^{n}|Y_i| + C_{F1}\bigg) > \varepsilon/2\bigg]$$

$$\leq \mathbf{P}\left[C_{F2}\zeta_\kappa r n^{-1}\sum_{i=1}^n(|Y_i|-\mathbf{E}|Y|)>\varepsilon/4\right]$$

$$\leq 2\exp(-a_1 n\varepsilon^2\zeta_\kappa^2)$$

for some constant $a_1>0$ and all sufficiently large $\kappa$ by Bernstein's inequality [Bosq (1998), page 22]. Also by Bernstein's inequality, there is a constant $a_2>0$ such that

$$\mathbf{P}[|\tilde{S}_{n\kappa 1}(\theta_{\kappa j})|>\varepsilon/2]\leq 2\exp(-a_2 n\varepsilon^2)$$

for all $n$, $\kappa$ and $j$. Therefore,

$$\begin{aligned}
\mathbf{P}_n &\leq 2[M\exp(-a_2 n\varepsilon^2)+\exp(-a_1 n\varepsilon^2)]\\
&\leq 2\exp\{-a_2 n\varepsilon^2\zeta_\kappa^2+2dC_\kappa n^\gamma[\log(C_\theta/r)+\tfrac{1}{2}\log(2C_\kappa d)+\tfrac{1}{2}\gamma\log n]\}\\
&\quad +2\exp(-a_1 n\varepsilon^2),
\end{aligned}$$

where $\gamma=4/15+\nu$. It follows that $\mathbf{P}_n\leq 4\exp(-an\varepsilon^2)$ for a suitable $a>0$ and all sufficiently large $n$. $\square$

Define

$$S_\kappa(\theta)=\mathbf{E}[S_{n\kappa}(\theta)]$$

and

$$\tilde{\theta}_\kappa=\underset{\theta\in\Theta_\kappa}{\arg\min}\,S_\kappa(\theta).$$

LEMMA 2. *For any $\eta>0$, $S_\kappa(\hat{\theta}_{n\kappa})-S_\kappa(\tilde{\theta}_\kappa)<\eta$ almost surely for all sufficiently large $n$.*

PROOF. For each $\kappa$, let $\mathcal{N}_\kappa\subset\mathbb{R}^{d(\kappa)}$ be an open set containing $\tilde{\theta}_\kappa$. Let $\bar{\mathcal{N}}_\kappa$ denote the complement of $\mathcal{N}_\kappa$ in $\Theta_\kappa$. Define $T_\kappa=\bar{\mathcal{N}}_\kappa\cap\Theta_\kappa$. Then $T_\kappa\subset\mathbb{R}^{d(\kappa)}$ is compact. Define

$$\eta=\min_{\theta\in T_\kappa}S_\kappa(\theta)-S_\kappa(\tilde{\theta}_\kappa).$$

Let $A_n$ be the event $|S_{n\kappa}(\theta)-S_\kappa(\theta)|<\eta/2$ for all $\theta\in\Theta_\kappa$. Then

$$A_n\Rightarrow S_\kappa(\hat{\theta}_{n\kappa})<S_{n\kappa}(\hat{\theta}_\kappa)+\eta/2$$

and

$$A_n\Rightarrow S_{n\kappa}(\tilde{\theta}_\kappa)<S_\kappa(\tilde{\theta}_\kappa)+\eta/2.$$

But $S_{n\kappa}(\hat{\theta}_{n\kappa})\leq S_{n\kappa}(\tilde{\theta}_\kappa)$ by definition, so

$$A_n\Rightarrow S_\kappa(\hat{\theta}_{n\kappa})<S_{n\kappa}(\tilde{\theta}_\kappa)+\eta/2.$$

Therefore,

$$A_n \Rightarrow S_\kappa(\hat{\theta}_{n\kappa}) < S_\kappa(\tilde{\theta}_\kappa) + \eta \Rightarrow S_\kappa(\hat{\theta}_{n\kappa}) - S_\kappa(\tilde{\theta}_\kappa) < \eta.$$

So $A_n \Rightarrow \hat{\theta}_{n\kappa} \in \mathcal{N}_\kappa$. Since $\mathcal{N}_\kappa$ is arbitrary, the result follows from Lemma 1 and Theorem 1.3.4 of Serfling [(1980), page 10].  □

Define

$$b_\kappa(x) = \mu + m(x) - P_\kappa(x)'\theta_\kappa$$

and

$$S_{\kappa 0}(\theta) = \mathbf{E}\{Y - F[P_\kappa(X)'\theta + b_\kappa(X)]\}^2.$$

Then

$$\theta_\kappa = \underset{\theta \in \Theta_\kappa}{\arg\min}\, S_{\kappa 0}(\theta).$$

LEMMA 3.  *For any* $\eta > 0$, $S_{\kappa 0}(\tilde{\theta}_\kappa) - S_{\kappa 0}(\theta_{\kappa 0}) < \eta$ *for all sufficiently large* $n$.

PROOF.  Observe that $|S_\kappa(\theta) - S_{\kappa 0}(\theta)| \to 0$ as $n \to \infty$ uniformly over $\theta \in \Theta_\kappa$ because $b_\kappa(x) \to 0$ for almost every $x \in \mathcal{X}$. For each $\kappa$, let $\mathcal{N}_\kappa \subset \mathbb{R}^{d(\kappa)}$ be an open set containing $\theta_{\kappa 0}$. Define $T_\kappa = \bar{\mathcal{N}}_\kappa \cap \Theta_\kappa$. Then $T_\kappa \subset \mathbb{R}^{d(\kappa)}$ is compact. Define

$$\eta = \min_{\theta \in T_\kappa} S_{\kappa 0}(\theta) - S_{\kappa 0}(\theta_{\kappa 0}).$$

By choosing a sufficiently small $\mathcal{N}_\kappa$, $\eta$ can be made arbitrarily small. Choose $n$ and, therefore, $\kappa$ large enough that $|S_\kappa(\theta) - S_{\kappa 0}(\theta)| < \eta/2$ for all $\theta \in \Theta$. Now proceed as in the proof of Lemma 2.  □

Define $Z_{\kappa i} = F'[\mu + m(X_i)]P_\kappa(X_i)$ and $\hat{Q}_\kappa = n^{-1}\sum_{i=1}^n Z_{\kappa i}Z_{\kappa i}'$. Then $Q_\kappa = \mathbf{E}\hat{Q}_\kappa$. Let $Z_{\kappa i}^k$ [$k = 1, \ldots, d(\kappa)$] denote the $k$th component of $Z_{\kappa i}$. Let $Z_\kappa$ denote the $n \times d(\kappa)$ matrix whose $(i, k)$ element is $Z_{\kappa i}^k$.

LEMMA 4.  $\|\hat{Q}_\kappa - Q_\kappa\|^2 = O_p(\kappa^2/n)$.

PROOF.  Let $Q_{ij}$ denote the $(i, j)$ element of $Q_\kappa$. Then

$$\mathbf{E}\|\hat{Q}_\kappa - Q_\kappa\|^2 = \sum_{k=1}^{d(\kappa)}\sum_{j=1}^{d(\kappa)}\mathbf{E}\left(n^{-1}\sum_{i=1}^n Z_{\kappa i}^k Z_{\kappa i}^j - Q_{kj}\right)^2$$

$$= \sum_{k=1}^{d(\kappa)}\sum_{j=1}^{d(\kappa)}\left(\mathbf{E}n^{-2}\sum_{i=1}^n\sum_{\ell=1}^n Z_{\kappa i}^k Z_{\kappa i}^j Z_{\kappa \ell}^k Z_{\kappa \ell}^j - Q_{kj}^2\right)$$

$$= \sum_{k=1}^{d(\kappa)} \sum_{j=1}^{d(\kappa)} \mathbf{E} n^{-2} \sum_{i=1}^{n} (Z_{\kappa i}^k)^2 (Z_{\kappa i}^j)^2 - n^{-1} \sum_{k=1}^{d(\kappa)} \sum_{j=1}^{d(\kappa)} Q_{kj}^2$$

$$\leq n^{-1} \mathbf{E} \left[ \sum_{k=1}^{d(\kappa)} (Z_{\kappa i}^k)^2 \sum_{j=1}^{d(\kappa)} (Z_{\kappa i}^j)^2 \right] = O(\kappa/n).$$

The lemma now follows from Markov's inequality. $\quad\square$

Define $\gamma_n = I(\lambda_{\kappa,\min} \geq c_\lambda/2)$, where $I$ is the indicator function. Let $\bar{U} = (U_1, \ldots, U_n)'$.

LEMMA 5. $\gamma_n \|\hat{Q}_\kappa^{-1} Z_\kappa' \bar{U}/n\| = O_p(\kappa^{1/2}/n^{1/2})$ as $n \to \infty$.

PROOF. For any $x \in \mathcal{X}$,

$$n^{-2} \mathbf{E}(\gamma_n \|\hat{Q}_\kappa^{-1/2} Z_\kappa' \bar{U}\|^2 | X = x) = n^{-2} \gamma_n \mathbf{E}(\bar{U}' Z_\kappa \hat{Q}_k^{-1} Z_\kappa \bar{U} | X = x)$$

$$= n^{-2} \mathbf{E}[\text{Trace}(Z_\kappa \hat{Q}_\kappa^{-1} Z_\kappa' \bar{U} \bar{U}') | X = x]$$

$$\leq n^{-2} \gamma_n C_V \text{Trace}(\hat{Q}_\kappa^{-1} Z_\kappa Z_\kappa')$$

$$= n^{-1} C_V \gamma_n d(\kappa) \leq C\kappa/n$$

for some constant $C < \infty$. Therefore, $\gamma_n \|\hat{Q}_\kappa^{-1/2} Z_\kappa' \bar{U}/n\| = O_p(\kappa^{1/2}/n^{1/2})$ by Markov's inequality. Now

$$\gamma_n \|\hat{Q}_\kappa^{-1} Z_\kappa' \bar{U}/n\| = \gamma_n [(\bar{U}' Z_\kappa/n) \hat{Q}_\kappa^{-1/2} \hat{Q}_\kappa^{-1} \hat{Q}_\kappa^{-1/2} (Z_\kappa' \bar{U}/n)]^{1/2}.$$

Define $\xi = \hat{Q}_\kappa^{-1/2} Z_\kappa' \bar{U}/n$. Let $\eta_1, \ldots, \eta_{d(\kappa)}$ and $q_1, \ldots, q_{d(\kappa)}$ denote the eigenvalues and eigenvectors of $\hat{Q}_\kappa^{-1}$. Let $\eta_{\max} = \max(\eta_1, \ldots, \eta_{d(\kappa)})$. The spectral decomposition of $\hat{Q}_\kappa^{-1}$ gives $\hat{Q}_\kappa^{-1} = \sum_{\ell=1}^{d(\kappa)} \eta_\ell q_\ell q_\ell'$, so

$$\gamma_n \|\hat{Q}_\kappa^{-1} Z_\kappa' \bar{U}/n\|^2 = \gamma_n \sum_{\ell=1}^{d(\kappa)} \eta_\ell \xi' q_\ell q_\ell' \xi$$

$$\leq \gamma_n \eta_{\max} \sum_{\ell=1}^{d(\kappa)} \xi' q_\ell q_\ell' \xi \leq \gamma_n \eta_{\max} \xi' \xi = O_p(\kappa/n). \quad\square$$

Define

$$B_n = \hat{Q}_\kappa^{-1} n^{-1} \sum_{i=1}^{n} F'[\mu + m(X_i)] Z_{\kappa i} b_{\kappa 0}(X_i).$$

LEMMA 6. $\|B_n\| = O(\kappa^{-2})$ with probability approaching 1 as $n \to \infty$.

PROOF. Let $\xi$ be the $n \times 1$ vector whose $i$th component is $F'[\mu + m(X_i)]b_{\kappa 0}(X_i)$. Then $B_n = \hat{Q}_\kappa^{-1} Z_\kappa' \xi / n$, and $\gamma_n \|B_n\|^2 = n^{-2} \gamma_n \xi' Z_\kappa \hat{Q}_\kappa^{-2} Z_\kappa' \xi$. Therefore, by the same arguments used to prove Lemma 5, $\gamma_n \|B_n\|^2 \leq C n^{-1} \gamma_n \xi' \xi = \gamma_n O(\kappa^{-4})$. The lemma follows from the fact that $\mathbf{P}(\gamma_n = 1) \to 1$ as $n \to \infty$. $\square$

PROOF OF THEOREM 1. To prove part (a), write

$$\begin{aligned}
S_{\kappa 0}(\hat{\theta}_{n\kappa}) - S_{\kappa 0}(\theta_\kappa) &= [S_{\kappa 0}(\hat{\theta}_{n\kappa}) - S_\kappa(\hat{\theta}_{n\kappa})] + [S_\kappa(\hat{\theta}_{n\kappa}) - S_\kappa(\tilde{\theta}_\kappa)] \\
&\quad + [S_\kappa(\tilde{\theta}_\kappa) - S_{\kappa 0}(\tilde{\theta}_\kappa)] + [S_{\kappa 0}(\tilde{\theta}_\kappa) - S_{\kappa 0}(\theta_\kappa)].
\end{aligned}$$

(7.1)

Given any $\eta > 0$, it follows from Lemmas 2 and 3 and uniform convergence of $S_\kappa$ to $S_{\kappa 0}$ that each term on the right-hand side of (7.1) is less than $\eta/4$ almost surely for all sufficiently large $n$. Therefore $S_{\kappa 0}(\hat{\theta}_{n\kappa}) - S_{\kappa 0}(\theta_\kappa) < \eta$ almost surely for all sufficiently large $n$. It follows that $\|\hat{\theta}_{n\kappa} - \theta_\kappa\| \to 0$ almost surely as $n \to \infty$ because $\theta_\kappa$ uniquely minimizes $S_\kappa$. Part (a) follows because uniqueness of the series representation of each function $m_j$ implies that $\|\theta_\kappa - \theta_{\kappa 0}\| \to 0$ as $n \to \infty$.

To prove the remaining parts of the theorem, observe that $\hat{\theta}_{n\kappa}$ satisfies the first-order condition $\partial S_{n\kappa}(\hat{\theta}_{n\kappa})/\partial\theta = 0$ almost surely for all sufficiently large $n$. Define $M_i = \mu + m(X_i)$ and $\Delta M_i = P_\kappa(X_i)'\hat{\theta}_{n\kappa} - M_i = P_\kappa(X_i)'(\hat{\theta}_{n\kappa} - \theta_{\kappa 0}) - b_{\kappa 0}(X_i)$. Then a Taylor series expansion yields

$$n^{-1} \sum_{i=1}^n Z_{\kappa i} U_i - (\hat{Q}_\kappa + R_{n1})(\hat{\theta}_{n\kappa} - \theta_{\kappa 0}) + n^{-1} \sum_{i=1}^n F'(M_i) Z_{\kappa i} b_{\kappa 0}(X_i) + R_{n2} = 0,$$

almost surely for all sufficiently large $n$. $R_{n1}$ is defined by

$$\begin{aligned}
R_{n1} = n^{-1} \sum_{i=1}^n \{ &-U_i F''(M_i) - U_i[F''(\tilde{\tilde{M}}_i) - F''(M_i)] \\
&+ [\tfrac{3}{2} F''(\tilde{M}_i) F'(M_i) + \tfrac{1}{2} F''(\tilde{M}_i) F''(\tilde{M}_i) \Delta M_i \\
&\qquad\qquad + \tfrac{1}{2} F''(\tilde{M}_i) F''(\tilde{M}_i)(\Delta M_i)^2] \Delta M_i \\
&- [F''(\tilde{\tilde{M}}_i) F'(M_i) - \tfrac{1}{2} F''(\tilde{M}_i) F'(\tilde{M}_i) \\
&\qquad\qquad + F''(\tilde{\tilde{M}}_i) F''(\tilde{M}_i) b_{\kappa 0}(X_i)] b_{\kappa 0}(X_i) \} \\
&\times P_\kappa(X_i) P_\kappa(X_i)',
\end{aligned}$$

where $\tilde{M}_i$ and $\tilde{\tilde{M}}_i$ are points between $P_\kappa(X_i)'\hat{\theta}_{n\kappa}$ and $M_i$. $R_{n2}$ is defined by

$$R_{n2} = -n^{-1} \sum_{i=1}^n \{ U_i F''(\tilde{\tilde{M}}_i) + U_i[F''(\tilde{\tilde{M}}_i) - F''(M_i)]$$

$$+ [F''(\tilde{\tilde{M}}_i)F'(M_i) - \tfrac{1}{2}F''(\tilde{M}_i)F'(M_i)]b_{\kappa 0}(X_i)$$

$$- \tfrac{1}{2}F''(\tilde{\tilde{M}}_i)F''(\tilde{M}_i)b_{\kappa 0}(X_i)^2\} P_\kappa(X_i)b_{\kappa 0}(X_i).$$

Now let $\xi$ denote either $\hat{Q}_\kappa^{-1}Z_\kappa'\bar{U}/n$ or $\hat{Q}_\kappa^{-1}[n^{-1}\sum_{i=1}^n F'(M_i)^2 P_\kappa(X_i) \times b_{\kappa 0}(X_i) + R_{n2}]$. Note that

$$\left\| n^{-1}\sum_{i=1}^n U_i F''(M_i)P_\kappa(X_i)P_\kappa(X_i)' \right\|^2 = O_p(\kappa^2/n).$$

Then

$$\gamma_n \|[(\hat{Q}_\kappa + R_{n1})^{-1} - \hat{Q}_\kappa^{-1}]\hat{Q}_\kappa \xi\|^2$$

$$= \gamma_n \|(\hat{Q}_\kappa + R_{n1})^{-1}R_{n1}\xi\|^2$$

$$= \operatorname{Trace}\{\gamma_n[\xi' R_{n1}(\hat{Q}_\kappa + R_{n1})^{-2}R_{n1}\xi]\}$$

$$= O_p(\|\xi' R_{n1}\|^2)$$

$$= O_p(\xi'\xi)O_p\Big\{\kappa^2/n + \int [P_\kappa(x)'(\hat{\theta}_{n\kappa} - \theta_{\kappa 0})]^2\,dx + \sup_{x \in \mathcal{X}} |b_{\kappa 0}(x)|^2\Big\}$$

$$= O_p(\xi'\xi)O_p(\kappa^2/n + \kappa\|\hat{\theta}_\kappa - \theta_{\kappa 0}\|^2 + \kappa^{-3}).$$

Setting $\xi = \hat{Q}_\kappa^{-1}Z_\kappa'\bar{U}/n$ and applying Lemma 5 yields $\|[(\hat{Q}_\kappa + R_{n1})^{-1} - \hat{Q}_\kappa^{-1}] \times Z_\kappa'\bar{U}/n\|^2 = O_p[\kappa^3/n + (\kappa^2/n)\|\hat{\theta}_{n\kappa} - \theta_{\kappa 0}\|^2 + 1/(n\kappa^2)]$. If $\xi = \hat{Q}_\kappa^{-1}[n^{-1} \times \sum_{i=1}^n F'(M_i)^2 P_\kappa(X_i)b_{\kappa 0}(X_i) + R_{n2}]$, then applying Lemma 6 and using the result $\|\hat{Q}_\kappa^{-1}R_{n2}\| = o_p(\kappa^{-2})$ yields

$$\left\| [(\hat{Q}_\kappa + R_{n1})^{-1} - \hat{Q}_\kappa^{-1}]\left[n^{-1}\sum_{i=1}^n F'(M_i)Z_{\kappa i}b_{\kappa 0}(X_i) + R_{n2}\right] \right\|^2$$

$$= O_p(\|\hat{\theta}_{n\kappa} - \theta_{\kappa 0}\|^2/\kappa + 1/\kappa^5).$$

It follows from these results that

$$\hat{\theta}_{n\kappa} - \theta_{\kappa 0} = n^{-1}\hat{Q}_\kappa^{-1}\sum_{i=1}^n F'[\mu + m(X_i)]P_\kappa(X_i)U_i$$

$$+ n^{-1}\hat{Q}_\kappa^{-1}\sum_{i=1}^n F'[\mu + m(X_i)]^2 P_\kappa(X_i)b_{\kappa 0}(X_i) + R_n,$$

where $\|R_n\| = O_p(\kappa^{3/2}/n + n^{-1/2})$. Part (d) of the theorem now follows from Lemma 4. Part (b) follows by applying Lemmas 5 and 6 to part (d). Part (c) follows from part (b) and Assumption A5(iii). $\square$

7.2. *Theorem* 2. This section begins with lemmas that are used to prove Theorem 2. For any $\tilde{x} \equiv (x^2, \ldots, x^d) \in [-1, 1]^{d-1}$, set $m_{-1}(\tilde{x}) = m_2(x^2) + \cdots + m_d(x^d)$, and $\bar{b}_{\kappa 0}(\tilde{x}) = \mu + m_{-1}(\tilde{x}) - \bar{P}_\kappa(\tilde{x})\bar{\theta}_{\kappa 0}$, where

$$\bar{P}_\kappa(x) = [1, 0, \ldots, 0, p_1(x^2), \ldots, p_\kappa(x^2), \ldots, p_1(x^d), \ldots, p_\kappa(x^d)]'$$

and

$$\bar{\theta}_{\kappa 0} = (\mu, 0, \ldots, 0, \theta_{21}, \ldots, \theta_{2\kappa}, \ldots, \theta_{d1}, \ldots, \theta_{d\kappa})'.$$

In other words, $\bar{P}$ and $\bar{\theta}_{\kappa 0}$ are obtained by replacing $p_1(x^1), \ldots, p_\kappa(x^d)$ with zeros in $P_\kappa$ and $\theta_{11}, \ldots, \theta_{1\kappa}$ with zeros in $\theta_{\kappa 0}$. Also define

$$\delta_{n1}(\tilde{x}) = n^{-1}\bar{P}_\kappa(\tilde{x})' Q_\kappa^{-1} \sum_{j=1}^n F'[\mu + m(X_j)] P_\kappa(X_j) U_j$$

and

$$\delta_{n2}(\tilde{x}) = n^{-1}\bar{P}_\kappa(\tilde{x})' Q_\kappa^{-1} \sum_{j=1}^n F'[\mu + m(X_j)]^2 P_\kappa(X_j) b_{\kappa 0}(X_j).$$

For $x^1 \in [-1, 1]$ and for $j = 0, 1$ define

$$H_{nj1}(x^1) = (nh)^{-1/2} \sum_{i=1}^n F'[\mu + m_1(x^1) + m_{-1}(\tilde{X}_i)]^2 (X_i^1 - x^1)^j$$
$$\times K_h(x^1 - X_i^1)\delta_{n1}(\tilde{X}_i),$$

$$H_{nj2}(x^1) = (nh)^{-1/2} \sum_{i=1}^n F'[\mu + m_1(x^1) + m_{-1}(\tilde{X}_i)]^2 (X_i^1 - x^1)^j$$
$$\times K_h(x^1 - X_i^1)\delta_{n2}(\tilde{X}_i)$$

and

$$H_{nj3}(x^1) = -(nh)^{-1/2} \sum_{i=1}^n F'[\mu + m_1(x^1) + m_{-1}(\tilde{X}_i)]^2 (X_i^1 - x^1)^j$$
$$\times K_h(x^1 - X_i^1)\bar{b}_{\kappa 0}(\tilde{X}_i).$$

Let $V(x) = \text{Var}(U|X = x)$.

LEMMA 7. *For $j = 0, 1$ and $k = 1, 2, 3$, $H_{njk}(x^1) = o_p(1)$ as $n \to \infty$ uniformly over $x^1 \in [-1, 1]$.*

PROOF. The proof is given only for $j = 0$. Similar arguments apply for $j = 1$. First consider $H_{n01}(x^1)$. We can write

$$H_{n01}(x^1) = \sum_{j=1}^n a_j(x^1) U_j,$$

where

$$a_j(x^1) = n^{-3/2} h^{-1/2} \sum_{i=1}^{n} F'[\mu + m_1(x^1) + m_{-1}(\tilde{X}_i)]^2$$
$$\times K_h(x^1 - X_i^1) \bar{P}_\kappa(\tilde{X}_i)' Q_\kappa^{-1} F'[\mu + m(X_j)] P_\kappa(X_j)$$
$$\equiv n^{-3/2} h^{-1/2} \sum_{i=1}^{n} K_h(x^1 - X_i^1) A_{ij}(x^1).$$

Define

$$\bar{a}(x^1) = \int F'[\mu + m_1(x^1) + m_{-1}(\tilde{x})]^2 \bar{P}_\kappa(\tilde{x}) f_X(x^1, \tilde{x}) \, d\tilde{x}.$$

Then arguments similar to those used to prove Lemma 1 show that $a_j(x^1) = (h/n)^{1/2}[\bar{a}(x^1) + r_n]' Q_\kappa^{-1} F'[\mu + m(X_j)] P_\kappa(X_j)$, where $r_n$ is uncorrelated with the $U_j$'s and $\|r_n\| = O[(\log n)/(nh)^{1/2}]$ uniformly over $x^1 \in [-1, 1]$ almost surely. Moreover, for each $x^1 \in [-1, 1]$, the components of $\bar{a}(x^1)$ are the Fourier coefficients of a function that is bounded uniformly over $x^1$. Therefore,

$$(7.2) \qquad \sup_{|x^1| \leq 1} \bar{a}(x^1)' \bar{a}(x^1) \leq M$$

for some finite constant $M$ and all $\kappa = 1, \ldots, \infty$. It follows from (7.2) and the Cauchy–Schwarz inequality that

$$\left\| \bar{a}(x^1)'(h/n)^{1/2} \sum_{j=1}^{n} U_j Q_\kappa^{-1} F'[\mu + m(X_j)] P_\kappa(X_j) \right\|^2$$
$$\leq M \left\| (h/n)^{1/2} \sum_{j=1}^{n} U_j Q_\kappa^{-1} F'[\mu + m(X_j)] P_\kappa(X_j) \right\|^2.$$

But

$$\mathbf{E} \left\| (h/n)^{1/2} \sum_{j=1}^{n} U_j Q_\kappa^{-1} F'[\mu + m(X_j)] P_\kappa(X_j) \right\|^2 = O(h),$$

so it follows from Markov's inequality that

$$\left\| \bar{a}(x^1)'(h/n)^{1/2} \sum_{j=1}^{n} U_j Q_\kappa^{-1} F'[\mu + m(X_j)] P_\kappa(X_j) \right\| = O_p(h^{1/2})$$

uniformly over $x^1 \in [-1, 1]$. This and $\|r_n\| = O[(\log n)/(nh)^{1/2}]$ establish the conclusion of the lemma for $j = 0$, $k = 1$.

We now prove the lemma for $j = 0$, $k = 2$. We can write

$$H_{n02}(x^1) = (nh)^{-1/2} \sum_{i=1}^{n} F'[\mu + m_1(x^1) + m_{-1}(\tilde{X}_i)]^2 K_h(x^1 - X_i^1) \bar{P}_\kappa(\tilde{X}_i)' B_n,$$

where

$$B_n = n^{-1} Q_\kappa^{-1} \sum_{j=1}^{n} F'[\mu + m(X_j)]^2 P_\kappa(X_j) b_{\kappa 0}(X_j).$$

Arguments like those used to prove Lemma 6 show that $\mathbf{E}\|B_n\|^2 = O(\kappa^{-4})$. Therefore,

$$\sup_{|x^1| \leq 1} |H_{n02}(x^1)| = \sup_{|x^1| \leq 1} (nh)^{-1/2} \sum_{i=1}^{n} K_h(x^1 - X_i^1) \cdot O_p(\kappa^{-2})$$

$$= O_p(n^{1/2} h^{1/2} \kappa^{-2})$$

$$= o_p(1).$$

For the proof with $j = 0$, $k = 3$, note that

$$\sup_{|x^1| \leq 1} |H_{n03}(x^1)| = \sup_{|x^1| \leq 1} (nh)^{-1/2} \sum_{i=1}^{n} K_h(x^1 - X_i^1) \cdot O_p(\kappa^{-2})$$

$$= o_p(1). \qquad \qquad \square$$

Let $P_\kappa(x^1, \tilde{X}) = [1, p_1(x^1), \ldots, p_\kappa(x^1), p_1(X^2), \ldots, p_\kappa(X^2), \ldots, p_1(X^d), \ldots, p_\kappa(X^d)]'$ and $b_\kappa(x^1, \tilde{X}_i) = \mu + m_1(x^1) + m_{-1}(\tilde{X}_i) - P_\kappa(x^1, \tilde{X}_i)\theta_{\kappa 0}$. Define

$$\delta_{n3}(x^1, \tilde{X}_i) = n^{-1} P_\kappa(x^1, \tilde{X}_i)' Q_\kappa^{-1} \sum_{j=1}^{n} F'[\mu + m(X_j)] P_\kappa(X_j) U_j$$

and

$$\delta_{n4}(x^1, \tilde{X}_i) = n^{-1} P_\kappa(x^1, \tilde{X}_i)' Q_\kappa^{-1} \sum_{j=1}^{n} F'[\mu + m(X_j)]^2 P_\kappa(X_j) b_{\kappa 0}(X_j).$$

Also, for $j = 0, 1$ define

$$L_{nj1}(x^1) = (nh)^{-1/2} \sum_{i=1}^{n} U_i F''[\mu + m_1(x^1) + m_{-1}(\tilde{X}_i)](X_i^1 - x^1)^j$$

$$\times K_h(x^1 - X_i^1)\delta_{n3}(x^1, \tilde{X}_i),$$

$$L_{nj2}(x^1) = (nh)^{-1/2} \sum_{i=1}^{n} U_i F''[\mu + m_1(x^1) + m_{-1}(\tilde{X}_i)](X_i^1 - x^1)^j$$

$$\times K_h(x^1 - X_i^1)\delta_{n4}(x^1, \tilde{X}_i),$$

$$L_{nj3}(x^1) = -(nh)^{-1/2} \sum_{i=1}^{n} U_i F''[\mu + m_1(x^1) + m_{-1}(\tilde{X}_i)](X_i^1 - x^1)^j$$
$$\times K_h(x^1 - X_i^1) b_{\kappa 0}(x^1, \tilde{X}_i),$$

$$L_{nj4}(x^1) = (nh)^{-1/2} \sum_{i=1}^{n} \{ F[\mu + m_1(X_i^1) + m_{-1}(X_i)]$$
$$- F[\mu + m_1(x^1) + m_{-1}(\tilde{X}_i)]\}$$
$$\times F''[\mu + m_1(x^1) + m_{-1}(\tilde{X}_i)](X_i^1 - x^1)^j$$
$$\times K_h(x^1 - X_i^1) \delta_{n3}(x^1, \tilde{X}_i),$$

$$L_{nj5}(x^1) = (nh)^{-1/2} \sum_{i=1}^{n} \{ F[\mu + m_1(X_i^1) + m_{-1}(X_i)]$$
$$- F[\mu + m_1(x^1) + m_{-1}(\tilde{X}_i)]\}$$
$$\times F''[\mu + m_1(x^1) + m_{-1}(\tilde{X}_i)](X_i^1 - x^1)^j$$
$$\times K_h(x^1 - X_i^1) \delta_{n4}(x^1, \tilde{X}_i)$$

and

$$L_{nj6}(x^1) = -(nh)^{-1/2} \sum_{i=1}^{n} \{ F[\mu + m_1(X_i^1) + m_{-1}(X_i)]$$
$$- F[\mu + m_1(x^1) + m_{-1}(\tilde{X}_i)]\}$$
$$\times F''[\mu + m_1(x^1) + m_{-1}(\tilde{X}_i)](X_i^1 - x^1)^j$$
$$\times K_h(x^1 - X_i^1) b_{\kappa 0}(x^1, \tilde{X}_i).$$

LEMMA 8. *As $n \to \infty$, $L_{njk}(x^1) = o_p(1)$ uniformly over $x^1 \in [-1, 1]$ for each $j = 0, 1$, $k = 1, \dots, 6$.*

PROOF. The proof is given only for $j = 0$. The arguments are similar for $j = 1$. By Theorem 1, $\delta_{n4}(x^1, \tilde{X}_i)$ is the asymptotic bias component of the stochastic expansion of $P_\kappa(x^1, \tilde{X}_i)(\hat{\theta}_{n\kappa} - \theta_{\kappa 0})$ and is $O_p(\kappa^{-3/2})$ uniformly over $(x^1, \tilde{X}_i) \in [-1, 1]^d$. This fact and standard bounds on

$$\sup_{|x^1| \leq 1} \sum_{i=1}^{n} |U_i| K_h(x^1 - X_i^1)$$

and

$$\sup_{|x^1| \leq 1} \sum_{i=1}^{n} K_h(x^1 - X_i^1)$$

establish the conclusion of the lemma for $j = 0$, $k = 2, 5$. For $j = 0$, $k = 3, 6$, proceed similarly using

$$\sup_{|x^1| \leq 1} |b_{\kappa 0}(x^1)| = O(\kappa^{-2}).$$

For $j = 0$, $k = 4$, one can use arguments similar to those made for $H_{n01}(x^1)$ in the proof of Lemma 7. It remains to consider $L_{n01}(x^1) = D_n(x^1)B_n$, where

$$D_n(x^1) = (nh)^{-1/2} \sum_{i=1}^{n} U_i F''[\mu + m_1(x^1) + m_{-1}(\tilde{X}_i)] K_h(x^1 - X_i^1) P_{\kappa}(x^1, \tilde{X}_i)$$

and

$$B_n = n^{-1} Q_{\kappa}^{-1} \sum_{j=1}^{n} F'[\mu + m(X_j)]^2 P_{\kappa}(X_j) U_j.$$

Now, $\mathbf{E}\|B_n\|^2 = O(\kappa n^{-1})$, and $D_n(x^1)$ contains elements of the form

$$(nh)^{-1/2} p_r(x^1) \sum_{i=1}^{n} U_i F''[\mu + m_1(x^1) + m_{-1}(\tilde{X}_i)] K_h(x^1 - X_i^1)$$

and

$$(nh)^{-1/2} \sum_{i=1}^{n} U_i p_r(X_i^\ell) F''[\mu + m_1(x^1) + m_{-1}(\tilde{X}_i)] K_h(x^1 - X_i^1)$$

for $0 \leq r \leq \kappa$, $2 \leq \ell \leq d$. These expressions can be bounded uniformly over $|x^1| \leq 1$ by terms that are $O_p[(\log n)^{1/2}|p_r(x^1)|]$ and $O_p[(\log n)^{1/2}]$, respectively. This gives

$$\sup_{|x^1| \leq 1} \|D_n(x^1)\|^2 = O_p(\kappa \log n).$$

Therefore,

$$\sup_{|x^1| \leq 1} |L_{n01}(x^1)|^2 \leq \sup_{|x^1| \leq 1} \|D_n(x^1)\|^2 \|B_n\|^2 = o_p(1). \qquad \square$$

LEMMA 9.  *The following hold uniformly over* $|x^1| \leq 1 - h$:

$$(nh)^{-1} S_{n01}''(x^1, \tilde{m}) = D_0(x^1) + o_p(1),$$
$$(nh)^{-1} S_{n21}''(x^1, \tilde{m}) = A_K h^2 D_0(x^1)[1 + o_p(1)]$$

*and*

$$(nh)^{-1} S_{n11}''(x^1, \tilde{m}) = h^2 A_K D_1(x^1)[1 + o_p(1)].$$

PROOF. This follows from Theorem 1(c) and standard bounds on

$$\sup_{|x^1| \leq 1} \sum_{i=1}^{n} U_i^r (X_i^1 - x^1)^s K_h(x^1 - X_i^1)$$

for $r = 0, 1$, $s = 0, 1, 2$. $\square$

Define $\Delta m_1(x^1) = \tilde{m}_1(x^1) - m_1(x^1)$, $\Delta m_{-1}(\tilde{x}) = \tilde{\mu} - \mu + \tilde{m}_{-1}(\tilde{x}) - m_{-1}(\tilde{x})$ and $\Delta m(x^1, \tilde{x}) = \Delta m_1(x^1) + \Delta m_{-1}(\tilde{x})$.

LEMMA 10. *The following hold uniformly over* $|x^1| \leq 1 - h$:

(a) $(nh)^{-1/2} S'_{n01}(x^1, \tilde{m}) = (nh)^{-1/2} S'_{n01}(x^1, m) + (nh)^{1/2} D_0(x^1) \Delta m_1(x^1) + o_p(1)$,

(b) $(nh)^{-1/2} S'_{n11}(x^1, \tilde{m}) = (nh)^{-1/2} S'_{n11}(x^1, m) + o_p(1)$.

PROOF. Only (a) is proved. The proof of (b) is similar. For each $i = 1, \ldots, n$, let $m^*(x^1, \tilde{X}_i)$ and $m^{**}(x^1, \tilde{X}_i)$ denote quantities that are between $\tilde{\mu} + \tilde{m}_1(x^1) + \tilde{m}_{-1}(\tilde{X}_i)$ and $\mu + m_1(x^1) + m_{-1}(\tilde{X}_i)$. The values of $m^*(x^1, \tilde{X}_i)$ and $m^{**}(x^1, \tilde{X}_i)$ may be different in different uses. A Taylor series expansion and Theorem 1(c) give

$$(nh)^{-1/2} S'_{n01}(x^1, \tilde{m}) = (nh)^{-1/2} S'_{n01}(x^1, m) + \sum_{j=1}^{4} J_{nj}(x^1)$$

$$+ n(nh)^{-1/2} O\left[ \sup_{(x^1, \tilde{x}) \in \mathcal{X}} \|\Delta m(x^1, \tilde{x})\|^3 \right]$$

$$= (nh)^{-1/2} S'_{n01}(x^1, m) + \sum_{j=1}^{4} J_{nj}(x^1) + o_p(1)$$

uniformly over $|x^1| \leq 1 - h$, where

$$J_{n1}(x^1) = 2(nh)^{-1/2} \sum_{i=1}^{n} F'[\mu + m_1(x^1) + m_{-1}(\tilde{X}_i)]^2 K_h(x^1 - X_i^1) \Delta m_1(x^1),$$

$$J_{n2}(x^1) = 2(nh)^{-1/2} \sum_{i=1}^{n} F'[\mu + m_1(x^1) + m_{-1}(\tilde{X}_i)]^2$$

$$\times K_h(x^1 - X_i^1) \Delta m_{-1}(\tilde{X}_i),$$

$$J_{n3}(x^1) = -2(nh)^{-1/2} \sum_{i=1}^{n} \{Y_i - F[\mu + m_1(x^1) + m_{-1}(\tilde{X}_i)]\}$$

$$\times F''[m^*(x^1, \tilde{X}_i)] K_h(x^1 - X_i^1) \Delta m(x^1, \tilde{X}_i),$$

$$J_{n4}(x^1) = 2(nh)^{-1/2} \sum_{i=1}^{n} F'[\mu + m_1(x^1) + m_{-1}(\tilde{X}_i)]$$
$$\times \{F''[m^*(x^1, \tilde{X}_i)] + 2F''[m^{**}(x^1, \tilde{X}_i)]\}$$
$$\times K_h(x^1 - X_i^1)[\Delta m(x^1, \tilde{X}_i)]^2.$$

It follows from Theorem 1(d) and Lemma 7 that $J_{n2}(x^1) = 2\sum_{k=1}^{3} H_{n0k}(x^1) + o_p(1) = o_p(1)$ uniformly over $|x^1| \leq 1 - h$. In addition, it follows from Theorem 1(c) that for some constant $C < \infty$,

$$J_{n4}(x^1) < C(nh)^{-1/2} \sum_{i=1}^{n} K_h(x^1 - X_i^1) \left[ \sup_{(x^1, \tilde{x}) \in \mathcal{X}} \|\Delta m(x^1, \tilde{x})\|^2 \right]$$
$$= O_p \left[ (nh)^{1/2} \sup_{(x^1, \tilde{x}) \in \mathcal{X}} \|\Delta m(x^1, \tilde{x})\|^2 \right] = o_p(1)$$

uniformly over $|x^1| \leq 1 - h$. Now consider $J_{n3}(x^1)$. It follows from Assumption A3(v) that

$$J_{n3}(x^1) = -2(nh)^{-1/2} \sum_{i=1}^{n} \{Y_i - F[\mu + m_1(x^1) + m_{-1}(\tilde{X}_i)]\}$$
$$\times F''[\mu + m_1(x^1) + m_{-1}(\tilde{X}_i)]$$
$$\times K_h(x^1 - X_i^1)\Delta m(x^1, \tilde{X}_i)$$
$$+ O_p \left[ (nh)^{1/2} \sup_{x \in \mathcal{X}} |\Delta m(x^1, \tilde{x})|^2 \right]$$
$$= \sum_{k=1}^{6} L_{n0k}(x^1) + O_p \left[ (nh)^{1/2} \sup_{x \in \mathcal{X}} |\Delta m(x^1, \tilde{x})|^2 \right] + o_p(1)$$

uniformly over $|x^1| \leq 1 - h$. Therefore, $J_{n3}(x^1) = o_p(1)$ uniformly by Lemma 8 and Theorem 1(c), and

$$(nh)^{-1/2} S'_{n01}(x^1, \tilde{m}) = (nh)^{-1/2} S'_{n01}(x^1, m) + J_{n1}(x^1) + o_p(1)$$

uniformly over $|x^1| \leq 1 - h$.

Now consider $J_{n1}(x^1)$. Set

$$\tilde{J}_{n1}(x^1) = 2(nh)^{-1/2} \sum_{i=1}^{n} F'[\mu + m_1(x^1) + m_{-1}(\tilde{X}_i)]^2 K_h(x^1 - X_i^1).$$

It follows from Theorem 2.37 of Pollard (1984) that $\tilde{J}_{n1}(x^1) - \mathbf{E}[\tilde{J}_{n1}(x^1)] = o(\log n)$ almost surely as $n \to \infty$. In addition, $\mathbf{E}[(nh)^{-1/2}\tilde{J}_{n1}(x^1)] = D(x^1) + O(h^2)$. Therefore,

$$J_{n1}(x^1) = (nh)^{1/2} D_0(x^1)\Delta m_1(x^1) + O[\log n \Delta m_1(x^1)]$$
$$= (nh)^{1/2} D_0(x^1)\Delta m_1(x^1) + o_p(1)$$

uniformly over $|x^1| \leq 1 - h$. $\quad\square$

PROOF OF THEOREM 2. By the definition of $\hat{m}_1(x^1)$,

$$
\begin{aligned}
&\hat{m}_1(x^1) - m_1(x^1) \\
(7.3) \quad &= \tilde{m}_1(x^1) - m_1(x^1) \\
&\quad - \frac{S''_{n21}(x^1, \tilde{m})S'_{n01}(x^1, \tilde{m}) - S''_{n11}(x^1, \tilde{m})S'_{n11}(x^1, \tilde{m})}{S''_{n01}(x^1, \tilde{m})S''_{n21}(x^1, \tilde{m}) - S''_{n11}(x^1, \tilde{m})^2}.
\end{aligned}
$$

Part (a) follows by applying Lemmas 9 and 10 to the right-hand side of (7.3). Define

$$
w = [nhD_0(x^1)]^{-1}\{-S'_{n01}(x^1, m) + [D_1(x^1)/D_0(x^1)]S'_{n11}(x^1, m)\}.
$$

Methods identical to those used to establish asymptotic normality of local-linear estimators show that $\mathbf{E}(n^{2/5}w) = \beta_1 + o(1)$, $\mathrm{Var}(n^{2/5}w) = V_1(x^1) + o(1)$ and $n^{2/5}[\hat{m}_1(x^1) - m_1(x^1)]$ is asymptotically normal, which proves part (b). $\quad\square$

PROOF OF THEOREM 4. It follows from Theorem 2(a) that

$$
n^{4/5} \int_{1-h \leq |x^1| \leq 1} w(x^1)[\hat{m}_1(x^1) - m_1(x^1)]^2 \, dx^1 = o_p(1).
$$

Now consider

$$
n^{4/5} \int_{|x^1| \leq 1-h} w(x^1)[\hat{m}_1(x^1) - m_1(x^1)]^2 \, dx^1.
$$

By replacing the integrand with the expansion of Theorem 2(a), one obtains a $U$-statistic in $U_i$ conditional on $X_1, \ldots, X_n$. This $U$-statistic has vanishing conditional variance. $\quad\square$

PROOF OF THEOREM 5. Use Theorem 2(a) to replace $\hat{m}_1$ with $m_1$ in the expression for $\hat{m}_1^{(\ell)}$. The result now follows from standard methods for bounding kernel estimators. $\quad\square$

DEPARTMENT OF ECONOMICS
NORTHWESTERN UNIVERSITY
EVANSTON, ILLINOIS 60208-2600
USA
E-MAIL: joel-horowitz@northwestern.edu

FAKULTÄT FÜR VOLKSWIRTSCHAFTSLEHRE
UNIVERSITÄT MANNHEIM
VERFÜGUNGSGEBÄUDE L7, 3-5
68131 MANNHEIM
GERMANY
E-MAIL: emammen@rumms.uni-mannheim.de